\documentclass{svjour3-ppt}
\smartqed
\journalname{}
\usepackage{graphicx,amssymb,amsmath}

\newcommand{\R}{\mathbb{R}}

\newcommand{\eps}{\varepsilon}
\newcommand{\be}[1]{\begin{equation}\label{#1}}
\newcommand{\ee}{\end{equation}}

\newcommand{\dx}{\,{\rm d}x}
\newcommand{\dy}{\,{\rm d}y}
\newcommand{\ds}{\,{\rm d}s}
\newcommand{\dt}{\,{\rm d}t}
\newcommand{\dr}{\,{\rm d}r}
\newcommand{\dz}{\,{\rm d}z}
\newcommand{\dxi}{\,{\rm d}\xi}
\newcommand{\dzeta}{\,{\rm d}\zeta}
\renewcommand{\(}{\left(}
\renewcommand{\)}{\right)}

\begin{document}
\title{Large mass self-similar solutions of the parabolic-parabolic Keller--Segel model of chemotaxis}

\titlerunning{Self-similar solutions of the Keller--Segel model}
\author{Piotr Biler \and Lucilla Corrias \and Jean Dolbeault}
\authorrunning{P. Biler, L. Corrias and J. Dolbeault}
\institute{
P. Biler \at Instytut Matematyczny, Uniwersytet Wroc\l awski, pl. Grunwaldzki 2/4, 50--384 Wroc\l aw, Poland. \email{Piotr.Biler@math.uni.wroc.pl}
\and
L. Corrias \at D\'epartement de Math\'ematiques, Universit\'e d'\'Evry Val d'Es\-sonne, rue du p\`ere Jarlan, \hbox{F-91025} \'Evry C\'edex, France. \email{lucilla.corrias@univ-evry.fr}
\and
J. Dolbeault \at Ceremade (UMR CNRS 7534), Universit\'e Paris-Dauphine, Place de Lattre de Tassigny, \hbox{F-75775} Paris C\'edex 16, France. \email{dolbeaul@ceremade.dauphine.fr}
}
\date{\today}
\maketitle

\begin{abstract} In two space dimensions, the parabolic-parabolic Keller--Segel system shares many properties with the parabolic-elliptic Keller--Segel system. In particular, solutions globally exist in both cases as long as their mass is less than $8\,\pi$. However, this threshold is not as clear in the parabolic-parabolic case as it is in the parabolic-elliptic case, in which solutions with mass above $8\,\pi$ always blow up.
Here we study forward self-similar solutions of the parabolic-parabolic Keller--Segel system and prove that, in some cases, such solutions globally exist even if their total mass is above $8\,\pi$, which is forbidden in the parabolic-elliptic case.

\keywords{Keller--Segel model \and chemotaxis \and self-similar solution \and nonlocal parabolic equations \and critical mass \and existence \and blowup}
\subclass{35B30 \and 35K40 \and 35K57 \and 35J60}
%
%
%
\end{abstract}

\section{Introduction}\label{Sec:intro}

The Keller--Segel model has been widely studied for almost forty years. It models the behavior of a slime mold of myxamoebae, {\em Dictyostelium Discoideum}, which have the peculiarity of organizing themselves to form aggregates by moving towards regions of a higher concentration of a chemoattractant. This chemoattractant, the {\em cyclic adenosine monophosphate}, is secreted by the amoebae themselves when they are lacking of nutrients. The Keller--Segel model is considered as a prototypical model for pattern formation in chemotaxis, and has attracted a lot of attention as a test case for more complex taxis phenomena driven by chemical substances. See \cite{MR2448428,MR1908133,MR2013508,MR2430318,MR2430317} for further references.

The simplest version of the model is made of two equations, one for the density of the amoebae and another one for the density of the chemoattractant. Both are parabolic, although an even simpler version has been widely considered by neglecting the time-dependence of the density of the chemoattractant. We shall refer to the complete version of the model as the parabolic-parabolic model, and to the latter as the parabolic-elliptic Keller--Segel model.

The precise dependence of the diffusion coefficients and of the chemosensitivity parameter depends on the context. All variants of the model involve diffusions in the equations for the density of the amoebae and for the density of the chemoattractant. The coupling is due to the fact that amoebae move according to the gradient of the chemoattractant, and that the emission of the chemoattractant is proportional to the density of amoebae. A crude insight into the main features of the model can be gained from the simplest case, that is when the nonlinear term in the equation is quadratic, but more realistic models should probably involve more complex nonlinearities.

Since the slime mold usually moves over a planar substrate, it makes sense to consider two-dimensional geometries. In some cases, boundary effects are important, but they are out of the purpose of this paper and we shall therefore assume that the model is set on the two-dimensional Euclidean plane.

In the parabolic-elliptic model, there is a critical mass, $8\,\pi$ after a proper adimensionalization, whose role is now rather well understood; see \cite{MR2103197,BDP}. Below such a~mass, the diffusion predominates, in the sense that amoebae are unable to emit enough chemoattractant to aggregate. On large times, the population diffuses and locally vanishes, although the behavior significantly differs from a pure diffusion. Above $8\,\pi$, at least one singularity appears in finite time, which is interpreted as the occurrence of an aggregate.

Since singularities are local, it is widely believed that $8\,\pi$ should also be a threshold between the diffusion dominated regime and the regime of aggregation also in the parabolic-parabolic model. This is certainly the case in some sense, for appropriate initial data, but the situation is not as simple as in the parabolic-elliptic case. It turns out that if, initially, the population of amoebae is scattered enough, and for a~well chosen initial distribution of the chemoattractant, there are solutions for which the diffusion predominates for large times, even for masses larger than $8\,\pi$. It is the purpose of this paper to establish such a fact, for a special class of solutions and in a certain range of the parameters of the model.
\medskip

In this paper, we consider the parabolic-parabolic Keller--Segel model
\begin{eqnarray}
\label{KSn}
&&n_t=\Delta n-\nabla\cdot(n\,\nabla c)\,,\\
\label{KSc}
&&\tau\,c_t=\Delta c+n\,,
\end{eqnarray}
for the densities $n$ and $c$ of, respectively, microorganisms ({\em e.g.} amoebae) and diffusing chemicals that they are secreting. Interesting mathematical questions are related to qualitative properties of problem \eqref{KSn}--\eqref{KSc} such as global in time existence {\em versus} finite time blowup of solutions describing chemotactic concentration phenomena. After the pioneering works of Keller and Segel, a huge literature has dealt with the mathematical modelling of chemotaxis and its analysis. We recommend the reading of \cite{MR2448428} for a recent review from both biological and mathematical points of view.

We shall consider the Keller--Segel system \eqref{KSn}--\eqref{KSc} for any $t>0$, $x\in\R^2$, supplemented with initial conditions $n_0$ and $c_0$. From now on we shall assume that $n_0$ and $c_0$ are nonnegative and that $n_0$ is integrable on $\R^2$. As a consequence, for solutions with sufficiently fast decay at infinity, the total mass is conserved, {\em i.e.},
\[
M:=\int_{\R^2}n(t,x)\dx=\int_{\R^2}n_0(x)\dx
\]
does not depend on $t$.

Throughout the paper, $\tau$ is a nonnegative parameter taking into account the difference of the time scales of the diffusive processes undergone by $n$ and $c$. The qualitative properties of $n$ and $c$ (such as the asymptotic behavior for large values of $t$) depend on $\tau$ and the stability of system \eqref{KSn}--\eqref{KSc} with respect to $\tau$ is expected, {\em i.e.} solutions of the parabolic-parabolic Keller--Segel system are expected to converge to those of parabolic-elliptic system when $\tau\searrow 0$. This has been recently proved, at least for solutions with a suitably small mass~$M$, in \cite{AR}. Here, we are interested in the differences between the parabolic-elliptic Keller--Segel system ($\tau=0$) and the parabolic-parabolic Keller--Segel system ($\tau>0$). Known results are briefly summarized as follows.

When $\tau=0$ in \eqref{KSc}, $M=8\,\pi$ is a threshold for existence {\em versus} blowup of the solution of \eqref{KSn}--\eqref{KSc}, see \cite{MR2103197,BDP,CC08}. Solutions globally exist for $M<8\,\pi$, while explosion in finite time may occur if $M>8\,\pi$. In the critical case $M=8\,\pi$, the solutions are known to be global in time but the density grows and mass concentration occurs in infinite time; see \cite{BKLN,BCM}.

For $\tau>0$, according to \cite{CC08}, solutions globally exist for any $M<8\,\pi$. However, it has not yet been proved that explosion occurs in finite time as soon as $M>8\,\pi$, for instance under some additional assumptions like a smallness condition on $\int_{\R^2}|x|^2\,n_0(x)\dx$. If $M=8\,\pi$, there is an infinite number of steady states (see \cite{BKLN}), but no other result is available, apart from self-similar solutions.
\medskip

Motivated by this lack of results for \eqref{KSn}--\eqref{KSc}, this paper deals with the existence of \emph{positive forward self-similar solutions} of \hbox{\eqref{KSn}--\eqref{KSc}}, {\em i.e.}, solutions which can be written~as
\begin{equation}
\label{eq:FSSS}
n(t,x)=\frac 1t\,u\(\frac x{\sqrt t}\)\quad\mbox{and}\quad c(t,x)=v\(\frac x{\sqrt t}\),
\end{equation}
with \emph{a large total mass} (that is, larger than $8\,\pi$). Indeed, since we are dealing with the two-dimensional case, any self-similar solution $n$ in $L^1(\R^2)$ preserves mass, {\em i.e.}, for each $t\ge 0$ 
%
\[
\int_{\R^2}n(t,x)\dx=\int_{\R^2}u(\xi)\dxi=M\,.
\]
Therefore, for any given $\tau>0$, we are interested in the optimal range of $M$ for the existence of such solutions, and in uniqueness or multiplicity issues for a given $M$ in the optimal range. Actually, our goal is double. The main one is to prove the above mentioned existence result. Second, we will give an as complete as possible review of the numerous existing results on the topic and also simplified, new proofs of them. For this reason, the remainder of the introduction will be primarily devoted to the state of the art on self-similar solutions.
\medskip

Self-similar solutions can be obtained through various approaches. The first method for the study of self-similar solutions (see for example \cite{B98} and the references therein) amounts to look for mild solutions of \eqref{KSn}--\eqref{KSc}, that is, solutions of
\begin{eqnarray*}
&&n(t,\cdot)={\rm e}^{(t-t_0)\Delta}\,n(t_0,\cdot)-\int_{t_0}^t\(\nabla {\rm e}^{(t-s)\Delta}\)\cdot\big(n(s,\cdot)\,\nabla c(s,\cdot)\big)\ds\,,\\
&&c(t,\cdot)={\rm e}^{\frac {t-t_0}\tau\Delta}\,c(t_0,\cdot)+\frac 1\tau\int_{t_0}^t{\rm e}^{\frac{t-s}\tau\Delta}\,n(s,\cdot)\ds\,,
\end{eqnarray*}
for any $t>t_0\ge 0$. Roughly speaking, such self-similar solutions are obtained by a~fixed point theorem. However, smallness conditions on the initial data are required in order to apply a contraction mapping principle; see \cite{N06}, where this method has been applied to \eqref{KSn}--\eqref{KSc} with $\tau=1$. Therefore, covering the whole range of masses for which solutions exists seems out of reach in this setting.
\medskip

Alternatively, one can prove the existence of self-similar solutions through the direct analysis of the elliptic system satisfied by $(u,v)$, {\em i.e.},
\begin{eqnarray}
\label{KSuSS}
&&\Delta u-\nabla\cdot\(u\,\nabla v-\frac12\,\xi\,u\)=0\,,\\
\label{KSvSS}
&&\Delta v+\frac\tau{2}\,\xi\cdot\nabla v+u=0\,,
\end{eqnarray}
where $\xi=x/\sqrt t$ and the differential operators in \eqref{KSuSS}--\eqref{KSvSS} are taken with respect to $\xi$. In this case, a natural functional space to be considered for both $u$ and $v$ is the subspace $C^2_0(\R^2)$ of functions in the space $C^2(\R^2)$ such that
\begin{equation*}
\lim_{|\xi|\to\infty}u(\xi)=0\quad\mbox{and}\quad\lim_{|\xi|\to\infty}v(\xi)=0\,.
\end{equation*}
For such classical solutions, equation \eqref{KSuSS} can be written equivalently as either
\[
\nabla\cdot\left[u\,\nabla\(\log u-v+\frac{|\xi|^2}4\)\right]=0\,,
\]
or
\[
\nabla\cdot\left[{\rm e}^v\,{\rm e}^{-|\xi|^2/4}\,\nabla \(u\,{\rm e}^{-v}\,{\rm e}^{|\xi|^2/4}\)\right]=0\,.
\]
Then, using the fact that $u$, $v$, and consequently $|\nabla v|$ are bounded, it has been proved in \cite{NSY02} that there exists a constant $\sigma$ such that
\begin{equation}
u(\xi)=\sigma\,{\rm e}^{v(\xi)}{\rm e}^{-\frac{|\xi|^2}{4}}
\label{u}
\end{equation}
for any $\xi\in\R^2$. Since $u$ is positive by the maximum principle, it follows that $\sigma$ is positive. As a consequence, $u\in L^1(\R^2)$, and the stationary system \eqref{KSuSS}--\eqref{KSvSS} reduces to a family of nonlinear elliptic equations for $v$, namely
\begin{equation}
\Delta v+\frac\tau{2}\,\xi\cdot\nabla v+\sigma\,{\rm e}^{v}\,{\rm e}^{-\frac{|\xi|^2}{4}}=0\,,
\label{eq:v}
\end{equation}
parametrized by $\sigma>0$. Again by the maximum principle applied to \eqref{eq:v}, the following upper bound for $v$ can be proved
\begin{equation}
v(\xi)\le C\,{\rm e}^{-\min\{1,\tau\}\frac{|\xi|^2}{4}},
\label{est:v}
\end{equation}
where $C$ is any positive constant such that $C\min\{1,\tau\}\ge\sigma\,{\rm e}^{\|v\|_\infty}$; see for instance~\cite{NSY02}. Therefore, $v\in L^1(\R^2)$ holds true for any solution of \eqref{eq:v} in $C^2_0(\R^2)$.

The range of $M$ for which self-similar solutions exist in $C^2_0(\R^2)$ gives an indication on the range of $M$ for which some solutions of \eqref{KSn}--\eqref{KSc} may globally exist. Self-similar solutions indeed provide explicit examples of global solutions, even with smooth initial data, up to a time-shift: take for instance $u$ and $v$ as the initial data for \eqref{KSn}--\eqref{KSc}. Moreover, if self-similar solutions describe the asymptotic behavior of any solution of \eqref{KSn}--\eqref{KSc} under appropriate conditions on initial data, then the ranges of global existence of solutions should be exactly the same. This property has been established in \cite{BDP} for $\tau=0$. In the case $\tau>0$, this might not be as simple as in the case $\tau=0$ if one can prove that blowup may occur for any $M>8\,\pi$. However, at least for initial data close enough to $u$ and~$v$, one can expect that the ranges of global existence are the same.

In view of our main goal, we are actually more interested in parametrizing the set of $C^2_0(\R^2)$ self-similar solutions in terms of mass rather than in terms of $\sigma$. This is possible using in \eqref{eq:v} the relation
\begin{equation}
M=\sigma\int_{\R^2}{\rm e}^{v(\xi)}\,{\rm e}^{-\frac{|\xi|^2}{4}}\dxi\,.
\label{sigma}
\end{equation}
However, by doing that, equation~\eqref{eq:v} becomes nonlocal, as was the original system \hbox{\eqref{KSn}--\eqref{KSc}}, and the problem is definitely more difficult to handle. Another not less important reason to consider a different but equivalent formulation of problem \eqref{KSuSS}--\eqref{KSvSS} is that the correspondence between $\sigma$ and $M$ is not clear due to the lack of uniqueness of solutions to \eqref{eq:v}; see Remark \ref{rk:sigma} at the end of Section \ref{Sec:Estimate}.

For the sake of completeness, we have to say that equation \eqref{eq:v}, written as
\[
\nabla\cdot\left({\rm e}^{\frac\tau4|\xi|^2}\,\nabla v\right)+\sigma\;{\rm e}^{v}\;{\rm e}^{\frac{\tau-1}4|\xi|^2}=0\,,
\]
has been studied using variational methods in \cite{MNY00,Y01}. The weighted functional space $H^1(\R^2;\exp(\frac\tau4|\xi|^2)\dxi)$ is then natural, but working in this space introduces a condition on the values of~$\tau$, which have to be in the interval $(0,2)$. Under such a restriction, it has been established that solutions exist if $0<\sigma<\sigma^*$, for some $\sigma^*>0$. These solutions are positive and belong to $C^2_0(\R^2)$, but due to the restriction on $\tau$, one has to look for alternative approaches.
\medskip

Another important and useful result has been obtained in \cite{NSY02} using the moving planes technique: any positive solution $v\in C^2_0(\R^2)$ of \eqref{eq:v} must be radially symmetric. As a consequence, system \eqref{KSuSS}--\eqref{KSvSS} reduces to the ODE system
\begin{eqnarray}
\label{KSur}
&&u'-u\,v'+\frac12\,r\,u=0\,,\\
\label{KSvr}
&&v''+\(\frac1r+\frac\tau2\,r\)\,v'+u=0\,,
\end{eqnarray}
where $u$ and $v$ are considered as functions of the radial variable $r=|\xi|$ only. Equations \eqref{u}--\eqref{eq:v} then become
\begin{eqnarray}
\label{KSur2}\nonumber
&&u(r)=\sigma\,{\rm e}^{v(r)}\,{\rm e}^{-r^2/4},\\
\label{KSvr2}
&&v''+\(\frac1r+\frac\tau2\,r\)v'+\sigma\,{\rm e}^{v}\,{\rm e}^{-r^2/4}=0\,.
\end{eqnarray}

Equation \eqref{KSvr2} has been studied in \cite{MMY99,NSY02}. More specifically, the authors proved in~\cite{MMY99} the existence of a positive decreasing solution of \eqref{KSvr2} endowed with the initial and integrability conditions
\begin{equation}
v'(0)=0\quad\mbox{and}\quad \int_0^{\infty}r\,v(r)\dr<\infty\,,
\label{cond1}
\end{equation}
for any $\tau>0$ and $\sigma>0$ such that $\sigma\,\frac{\log\tau}{\tau-1}<1/{\rm e}$ (see Remark \ref{rk:CS}). However, such a condition does not determine the optimal range neither for the parameter $\sigma$ nor for $M$.

It is worth noticing that the boundary conditions \eqref{cond1} and the following ones,
\begin{equation}
v'(0)=0\quad\mbox{and}\quad \lim_{r\to\infty}v(r)=0\,,
\label{cond2}
\end{equation}
are equivalent for classical decreasing solutions. Indeed, \eqref{cond1} implies \eqref{cond2} and the converse holds true by~\eqref{est:v}. Using \eqref{cond2}, equation \eqref{KSvr2} turns out to be equivalent to
\begin{eqnarray}
\label{KSwr}
&&w''+\(\frac1r+\frac\tau2\,r\)w'+ {\rm e}^{w}\,{\rm e}^{-r^2/4}=0\,,\\
\label{cond3}
&&w'(0)=0\quad\mbox{and}\quad w(0)=s\,,
\end{eqnarray}
for some shooting parameter $s\in\R$. Indeed, if $w(r;s)$ is a classical solution of \eqref{KSwr}--\eqref{cond3} for a given $s\in\R$, then $w(\infty;s)=\lim_{r\to\infty}w(r;s)$ exists and is finite and $v(r)=w(r;s)-w(\infty;s)$ is a classical solution of \eqref{KSvr2}--\eqref{cond2} with $\sigma={\rm e}^{w(\infty;s)}$. Conversely, if $v$ is a classical solution of \eqref{KSvr2}--\eqref{cond2}, then $w(r;s)=v(r)+\log\sigma$ is a classical solution of \eqref{KSwr}--\eqref{cond3} with $s=v(0)+\log\sigma$ and again $\sigma={\rm e}^{w(\infty;s)}$ holds true. It follows that all solutions of \eqref{KSvr2}--\eqref{cond2} can be parametrized in terms of~$s$. See \cite{NSY02} for more details. Using this equivalence, the authors of \cite{NSY02} analyze the structure of the set of solutions of \eqref{KSvr2}--\eqref{cond2} seen as a one-parameter family; see Remark \ref{rk:sigma} at the end of Section \ref{Sec:Estimate}. Computations presented in Figs. 1 have been based on this parametrization of the solution set.

Last but not least, the parametrization of the solutions of \eqref{KSwr}--\eqref{cond3} in terms of $s$ allows us to parametrize the total mass $M$ in term of $s$ by
\begin{equation}
M(s)= 2\,\pi\int_0^\infty {\rm e}^{w(r;s)}\,{\rm e}^{-r^2/4}\,r\dr\,.
\label{eq:M(s)}
\end{equation}
But again, this does not provide an explicit computation for the optimal range of $M$. Computations presented in Fig. 2 (left) have also been based on this parametrization of $M$.
\medskip

Being this the state of the art, we will establish that the formulation of system \hbox{\eqref{KSur}--\eqref{KSvr}} in terms of \emph{cumulated densities} is better adapted to the qualitative description of $u$ and $v$. This is a classical technique used previously, for example, in the context of the parabolic-elliptic Keller--Segel system and astrophysical models; see \cite{B98,BKLN} and further references therein. For $\tau>0$, many qualitative properties of the solutions can still be proved in this framework. These will allow us to build \emph{positive forward self-similar solutions} of \hbox{\eqref{KSn}--\eqref{KSc}} satisfying \eqref{eq:FSSS}, which have an arbitrarily large mass when $\tau$ is large enough. The obtained results are summarized in the theorem below. One may interpret it by saying that the diffusion of $c$ described by \eqref{KSc} for positive large $\tau$ and some $M>8\,\pi$ may prevent the blowup of the solutions of the parabolic-parabolic Keller--Segel system. This is a major difference with the parabolic-elliptic case $\tau=0$, for which the response of $c$ to the variations of $n$ being instantaneous, any smooth solution with mass $M>8\,\pi$ must concentrate and blow~up in finite time.

\begin{theorem}
\label{Thm:Simplified}
For any $M>0$, there exists some $\tilde\tau(M)\ge 0$ such that for any $\tau\ge\tilde\tau(M)$ there is at least one solution $(u,v)$ of \eqref{KSuSS}--\eqref{KSvSS} in $(C^2_0(\R^2))^2$ with $u>0$ of mass $M$ and $v>0$. If $M<8\,\pi$, $\tilde\tau(M)=0$. If $M>8\,\pi$, $\tilde\tau(M)$ is positive and there are at least two solutions, except for the maximal possible value of~$M$. All solutions are radial, nonincreasing, with fast decay at infinity, and hence attain their maximum at $x=0$. They are uniquely determined by $a:=u(0)/2$, which in turn uniquely determines $M=M(a,\tau)$. Moreover, $\lim_{a\to\infty}M(a,\tau)=8\,\pi$, while, as $a\to\infty$, the corresponding solution $u$ concentrates into a Dirac delta distribution, up to the factor $8\,\pi$, and $v(0)=\|v\|_{L^\infty(\R^2)}$ becomes arbitrarily large.
\end{theorem}

This paper is organized as follows. We shall first establish the main {\em a priori} estimates for Theorem~\ref{Thm:Simplified} in the next section. The framework of cumulated densities is developed in Section~\ref{Sec:Cumulated}, which also contains more detailed statements than the ones of Theorem~\ref{Thm:Simplified}. The remaining {\em a priori} estimates and proofs are given in Sections \ref{Sec:qualitative} and \ref{Sec:Proofs}, respectively. Section \ref{Sec:Numerics} is devoted to some numerical results and Section \ref{Sec:Conclusions} to concluding remarks.

\section{Large mass positive forward self-similar solutions}
\label{Sec:Estimate}

Before restating the question of self-similar solutions in terms of cumulated densities, let us establish the key {\em a priori} estimate for Theorem~\ref{Thm:Simplified}, which proves that these solutions may have an arbitrary large mass when $\tau$ is large enough. This result is entirely new. Such an estimate can be obtained both from equation \eqref{KSvr2} and from the cumulated densities formulation. In this section, we shall establish this {\em a priori} estimate in the first setting. It will be translated in the cumulated densities framework in Section~\ref{Sec:qualitative}.

From now on, we shall parametrize $M$ in term of $a$ and $\tau$, {\em i.e.} $M=M(a,\tau)$, where $a=u(0)/2$ will be the shooting parameter in the cumulated densities shooting problem, see \eqref{phi2}--\eqref{S} and \eqref{IC}--\eqref{infcond} below.

A positive classical solution $v$ of \eqref{KSvr2}, \eqref{cond2} solves
\[
\(r\,{\rm e}^{\tau\,r^2/4}\,v'\)'+\sigma\,r\,{\rm e}^{(\tau-1)\,r^2/4}\,{\rm e}^v=0\,,
\]
which, after an integration on $(0,r)$, gives
\begin{equation}
v'(r)=-\frac\sigma{r}\,{\rm e}^{-\tau\,r^2/4}\int_0^r{\rm e}^{(\tau-1)\,z^2/4}\,{\rm e}^{v(z)}\,z\dz\,.
\label{eq:v'}
\end{equation}
As a consequence, $v'$ is nonpositive, so that $v(z)\le v(0)$ for any $z\ge 0$ and, for $\tau\ne1$,
\begin{equation}
v'(r)\ge-\frac\sigma{r}\,{\rm e}^{-\tau\,r^2/4}\,{\rm e}^{v(0)}\int_0^r{\rm e}^{(\tau-1)\,z^2/4}\,z\dz=
-\frac{2}{\tau-1}\,\frac\sigma r\,{\rm e}^{v(0)}\left({\rm e}^{-r^2/4}-{\rm e}^{-\tau\,r^2/4}\right)\,.
\label{est:v'}
\end{equation}
We observe that
\[
\frac{\rm d}{{\rm d}\tau}\int_0^\infty\left({\rm e}^{-r^2/4}-{\rm e}^{-\tau\,r^2/4}\right)\frac{2\,{\rm d}r}r=\int_0^\infty{\rm e}^{-\tau\,r^2/4}\,\frac r2\,{\rm d}r=\frac 1\tau\,.
\]
Hence, after one more integration of \eqref{est:v'} on $(0,\infty)$, we get, for any $\tau\ne1$,
\begin{equation}
v(0)\le\sigma\,{\rm e}^{v(0)}\,I(\tau)\quad\mbox{with}\quad I(\tau):=\frac{\log\tau}{\tau-1}\,.
\label{est:v(0)}
\end{equation}
Actually, it is easy to check that estimate \eqref{est:v(0)} holds true also for $\tau=1$ with $I(1)=1$. Since from \eqref{u} we have
\begin{equation}
\sigma\,{\rm e}^{v(0)}=u(0)=2\,a\,,
\label{est:sigma-v(0)}
\end{equation}
it has been proved that for each $\tau>0$,
\begin{equation}
0=\lim_{z\to\infty}v(z)\le v(r)\le v(0)\le2\,a\,I(\tau)
\label{est2:v(0)}
\end{equation}
for any $r\in\R_+$. On the other hand, by \eqref{sigma}, \eqref{cond2} and \eqref{est2:v(0)}, mass can be estimated for any positive $a$ and $\tau$~by
\begin{equation}
M=2\,\pi\,\sigma\int_0^\infty {\rm e}^{v(r)}\,{\rm e}^{-r^2/4}\,r\dr\ge2\,\pi\,\sigma\int_0^\infty {\rm e}^{-r^2/4}\,r\dr=4\,\pi\,\sigma
\ge8\,\pi\,a\,{\rm e}^{-2\,a\,I(\tau)}
\label{est:Mgrand}
\end{equation}
using \eqref{est:sigma-v(0)}. As a function of $a$, $\widetilde M(a,\tau):=8\,\pi\,a\,{\rm e}^{-2\,a\,I(\tau)}$ achieves its maximum at $a_*(\tau):=\frac 1{2\,I(\tau)}$, which proves that $M=M(a,\tau)$ verifies for each $\tau>0$
\begin{equation*}
\max_{a>0}M(a,\tau)\ge\widetilde M(a_*(\tau),\tau)=\frac{4\,\pi}{{\rm e}\,I(\tau)}\:,
\label{est2:Mgrand}
\end{equation*}
and it is clear that the right hand side can be made arbitrarily large for $\tau$ large enough. Hence, the corresponding density $u(r)=\sigma\,{\rm e}^{v(r)}{\rm e}^{-r^2/4}$ has mass $M>8\,\pi$ if $\frac{4\,\pi}{{\rm e}\,I(\tau)}>8\,\pi$, that is for any $\tau>\bar\tau$ with $\bar\tau$ such that $I(\bar\tau)=\frac1{2\,{\rm e}}$, {\em i.e.} $\bar\tau\approx 16.1109$. Also observe that for any $\tau>\bar\tau$ the density $u$ corresponding to $a=a_*(\tau)$ satisfies $u(0)=2\,a_*(\tau)>2\,{\rm e}$. Finally, using $v(z)\ge v(r)$ in \eqref{eq:v'} and integrating the inequality on $(0,\infty)$, one obtains ${\rm e}^{-v(0)}-\lim_{r\to\infty}{\rm e}^{-v(r)}\le -\,\sigma\,I(\tau)$, for any $\tau>0$. As a consequence, using \eqref{est:sigma-v(0)} and $\lim_{r\to\infty}{\rm e}^{-v(r)}=1$, we obtain that
\[
1-{\rm e}^{v(0)}\le -\,\sigma\,I(\tau)\,{\rm e}^{v(0)}=-\,2\,a\,I(\tau)\,.
\]
This gives the estimate
\begin{equation}
v(0)>\log(2\,a\,I(\tau)+1)\,,
\label{est3:v(0)}
\end{equation}
which implies that $v(0)$ becomes arbitrarily large as $a\to\infty$, for any $\tau>0$.

Estimates \eqref{est:Mgrand} and \eqref{est3:v(0)} can be read also as lower and upper bounds for $\sigma=2\,a\,{\rm e}^{-v(0)}$, namely
\begin{equation}
2\,a\,{\rm e}^{-2\,a\,I(\tau)}\le\sigma\le\min\left\{\frac M{4\,\pi}\,,\,\frac{2\,a}{2\,a\,I(\tau)+1}\right\}\,,
\label{est:sigma}
\end{equation}
hence showing that $\sigma$ takes arbitrarily large values for $\tau$ large enough.

\begin{remark}
\label{rk:sigma} \
Estimates \eqref{est:sigma} on $\sigma$ are new. The authors of \cite{NSY02} analyzed the map $s \mapsto\sigma(s)$, where $s$ is the shooting parameter defined in \eqref{cond3}, and they proved that it is a continuous map from $\R$ into $\R_+$ with $\lim_{s\to\pm\infty}\sigma(s)=0$. Therefore, $\sigma$ must be bounded for any fixed $\tau$ by $\sigma^*=\sigma(s^*)$, for some $s^*\in\R$, and problem \eqref{KSvr2}--\eqref{cond2} admits no solution for $\sigma>\sigma^*$, at least one solution for $\sigma=\sigma^*$ and finally (at least) two distinct solutions for $0<\sigma<\sigma^*$. However, estimates on $\sigma$ (or $\sigma^*$) were missing.
\end{remark} 
\begin{remark}
\label{rk:CS} \
Estimate \eqref{est:v(0)} says that, for any fixed $\sigma>0$ and $\tau>0$, $v(0)$ satisfies
\[
v(0)-\,\sigma\,I(\tau)\,{\rm e}^{v(0)}\le0\,.
\]
Since the function $x\, \mapsto\, x-\sigma\,I(\tau)\,{\rm e}^x$ is strictly concave and attains the maximum in $x=-\log(\sigma\,I(\tau))$, we deduce that whenever $\sigma\,I(\tau)<1/{\rm e}$, there exists an open interval $J\subset\R_+$ of non existence of solutions of \eqref{KSvr2} satisfying \eqref{cond2}, with $v(0)\in J$. On the other hand if $\sigma\,I(\tau)\ge 1/{\rm e}$, the above inequality induces no restriction on $v(0)$.
\end{remark}

\section{Cumulated densities and main results}
\label{Sec:Cumulated}

Let us introduce the \emph{cumulated densities} formulation of the parabolic-parabolic Keller--Segel model as in \cite{B98}, in terms of the functions $u$ and $v$ which solve problem \eqref{KSur}--\eqref{KSvr}, by defining
\begin{eqnarray*}
&&\phi(y):=\frac 1{2\,\pi}\int_{B(0,\sqrt y)}u(\xi)\dxi=\int_0^{\sqrt y}r\,u(r)\dr\,,\\
&&\psi(y):=\frac 1{2\,\pi}\int_{B(0,\sqrt y)}v(\xi)\dxi=\int_0^{\sqrt y}r\,v(r)\dr\,.
\end{eqnarray*}
Using the relations
\begin{eqnarray}
\label{Derphi}
&&\phi'(y)=\frac12\,u\(\sqrt y\)\quad\mbox{and}\quad\phi''(y)=\frac 1{4\,\sqrt y}\,u'\(\sqrt y\)\,,\\
\label{Derpsi}\nonumber
&&\psi'(y)=\frac12\,v\(\sqrt y\)\quad\mbox{and}\quad\psi''(y)=\frac 1{4\,\sqrt y}\,v'\(\sqrt y\)\,,
\end{eqnarray}
it follows from \eqref{KSur}--\eqref{KSvr} that the cumulated densities $\phi$ and $\psi$ solve the second order ODE system
\begin{eqnarray}
\label{phi1}
&&\phi''+\frac14\,\phi'-2\,\phi'\psi''=0\,,\\
\label{psi}
&&4\,y\,\psi''+\tau\,y\,\psi'-\tau\,\psi+\phi=0\,,
\end{eqnarray}
where \eqref{KSvr} has been multiplied by $r$ and integrated on $(0,\sqrt y)$. Observing that equation~\eqref{psi} can be written as
\[
4\,(y\,\psi'-\psi)'+\tau\,(y\,\psi'-\psi)+\phi=0\,,
\]
and defining $S(y):=4\,(\psi(y)-y\,\psi'(y))'=-4\,y\,\psi''(y)=-{\sqrt y}\,v'(\sqrt y)$ as in \cite{B06,NSY02}, system \eqref{phi1}--\eqref{psi} becomes, after a differentiation of \eqref{psi} with respect to $y$, a first order system in the $(\phi',S)$ variables
\begin{eqnarray}
\label{phi2}
&&\phi''+\frac14\,\phi'+\frac1{2\,y}\,\phi'S=0\,,\\
\label{S}
&&S'+\frac\tau4\,S=\phi'\,.
\end{eqnarray}
The last formulation of the ODE system can be equivalently written as a single integro-differential equation, hence nonlocal, for $\phi'$,
\begin{equation}
\phi''+\frac14\,\phi'+\frac1{2\,y}\,\phi'\,{\rm e}^{-\tau\,y/4}\int_0^y{\rm e}^{\tau\,z/4}\,\phi'(z)\dz=0\,,
\label{phi3}
\end{equation}
since, by \eqref{S},
\begin{equation}
S(y)={\rm e}^{-\tau\,y/4}\int_0^y{\rm e}^{\tau\,z/4}\,\phi'(z)\dz\,,
\label{S1}
\end{equation}
and as a single, local but nonlinear second order ODE for $S$,
\begin{equation*}
S''+\frac14\,(\tau+1)\,S'+\frac\tau{16}\,S+\frac1{2\,y}\(S\,S'+\frac\tau4\,S^2\)=0\,,
\end{equation*}
which is obtained by differentiating \eqref{S}. We will use in the sequel all these formulations in order to get {\em a priori} estimates.

For any positive self-similar solution $(u,v)\in (C^2_0(\R^2))^2$, the natural initial conditions for \eqref{phi2}--\eqref{S} are
\begin{equation}
\phi(0)=0\,,\quad\phi'(0)=a>0\quad\mbox{and}\quad S(0)=0\,,
\label{IC}
\end{equation}
in view of the definition of $\phi$ and of \eqref{S1}. Moreover, for any self-similar solution $u\in L^1(\R^2)$, the corresponding cumulated density $\phi$ satisfies the boundary condition
\begin{equation}
\phi(\infty):=\lim_{y\to\infty}\phi(y)=\frac {M(a,\tau)}{2\,\pi}\,.
\label{infcond}
\end{equation}

The problem is now formulated in terms of a shooting parameter problem \eqref{phi2}--\eqref{S}, \eqref{IC}, with a new shooting parameter $a$ which is directly related to the concentration of the self-similar density $u$ around the origin, since $a=u(0)/2$. This has been obtained in Section~\ref{Sec:Estimate} and will be made more precise below. Let us observe that the relation between $a$ and the shooting parameter $s$ defined in \eqref{cond3} is $2\,a={\rm e}^s$, since $s=v(0)+\log\sigma$. Thus, a one-to-one relation is established between the initial valued problems \eqref{phi2}--\eqref{S}, \eqref{IC} and \eqref{KSwr}--\eqref{cond3} as soon as an existence and uniqueness result is established for one of them. Moreover, we have
\[
v(0)=v(\sqrt y)+\frac12\int_0^{y}\frac {S(z)}z\dz\,,
\]
and the boundary condition $\lim_{r\to\infty}v(r)=0$ is equivalent to
\begin{equation}\label{eq:v(0)}
v(0)=\frac12\int_0^{\infty}\frac {S(z)}z\dz\,.
\end{equation}
We also have: $\sigma=\lim_{r\to\infty}u(r)\,{\rm e}^{{r^2}/{4}}=2\lim_{y\to\infty}\phi'(y)\;{\rm e}^{\,y/4}$. Hence we can repara\-me\-trize $v(0)$ and $\sigma$ in terms of $a$ and $\tau$.
\medskip

The main statements we are going to prove are summarized in the following theorems. The {\em a priori} estimates will be established in Section \ref{Sec:qualitative}. The proofs will be given in Section~\ref{Sec:Proofs}. We shall say that $(\phi,S)$ is a {\em positive} solution if both $\phi$ and $S$ are positive functions.

\begin{theorem} For any $(a,\tau)\in\R^2_+$ there exists a unique positive solution $(\phi,S)$ of \eqref{phi2}--\eqref{S}, \eqref{IC} such that $\phi\in C^2(0,\infty)\cap C^1[0,\infty)$ and $S\in C^1[0,\infty)$. Moreover, for any fixed $\tau>0$, $\phi\in C^2[0,\infty)$, the maps $a\in\R_+\mapsto (\phi,S)$ and $a\in\R_+\mapsto M(a,\tau)\in\R_+$ are continuous and
\[
g(a,\tau)\le \frac {M(a,\tau)}{2\,\pi}\le f(a,\tau)\,,
\]
where
\begin{equation}
f(a,\tau)=\left\{
\begin{array}{lll}
\min\{4,4\,a\}\quad&\mbox{if}\quad&\tau\in\left(0,\frac12\right]\,,\\ \\
\min\left\{4\,a,\frac23\,\pi^2\right\}\quad&\mbox{if}\quad&\tau\in\left(\frac12,1\right]\,,\\ \\
\min\left\{4\,a,\frac23\,\pi^2\,\tau,4\,(\tau+1)\right\}\quad&\mbox{if}\quad&\tau>1\,,
\end{array}
\right.
\label{f}
\end{equation}
and
\begin{equation}
g(a,\tau)=\left\{
\begin{array}{rcl}
\max\left\{4\,a\;{\rm e}^{-2\,a\,\frac{\log\tau}{\tau-1}},\frac{4\,a\,\tau}{a+\tau}\right\}\quad&\mbox{if}\quad&\tau\in(0,1]\,,\\ \\
\max\left\{4\,a\,{\rm e}^{-2\,a\,\frac{\log\tau}{\tau-1}},\frac{4\,a}{a+1}\right\}\quad&\mbox{if}\quad&\tau>1\,.
\end{array}
\right.
\label{g}
\end{equation}
\label{the:existence}
\end{theorem}

For consistency, it is worth noticing that the inequality $g(a,\tau)\le f(a,\tau)$ holds for all $\tau>0$ and $a>0$.

\begin{theorem}
Given any fixed $\tau>0$, for any positive sequence $\{a_k\}$ such that $a_k\to\infty$ as $k\to\infty$, there exists a sequence of positive self-similar solutions $(u_k,v_k)\in (C^2_0(\R^2))^2$ satisfying \eqref{KSuSS}--\eqref{KSvSS} and $u_k(0)=2\,a_k$, $v_k'(0)=0$ such that
\[
u_k\rightharpoonup 8\,\pi\,\delta_0\quad\mbox{as}\quad k\to\infty\;
\]
in the sense of weak convergence of measures.
Moreover, $\lim_{k\to\infty}\int_{\R^2}u_k\dx=8\,\pi$ and $\lim_{k\to\infty}\|v_k\|_{L^\infty(\R^2)}~=~\infty$.
\label{the:BU}
\end{theorem}

Theorem~\ref{the:BU} has already been proved in \cite[Th.~2, (iii)]{NSY02} using a classical result by Brezis and Merle in \cite{BM}. However, here we shall give a simplified and quite direct proof using the cumulated densities formulation.

\begin{theorem}
For any fixed $\tau>0$ there exists $M^*=M^*(\tau)\ge 8\,\pi$ such that problem \eqref{phi2}--\eqref{S} with the boundary conditions
\begin{equation*}
\phi(0)=0\,,\quad\lim_{y\to\infty}\phi(y)=\frac M{2\,\pi}\,,\quad S(0)=0\,,
\label{IC2}
\end{equation*}
has no positive solution $(\phi,S)\in C^2[0,\infty)\times C^1[0,\infty)$ if $M>M^*$ and has at least one positive solution $(\phi,S)\in C^2[0,\infty)\times C^1[0,\infty)$ in the following cases:
\begin{description}
\item[(i)] $M\in(0,M^*]$ if $M^*>8\,\pi$,
\item[(ii)] $M\in(0,M^*)$ if $M^*=8\,\pi$.
\end{description}
Moreover, there exist $1/2<\tau^*\le \tau^{**}$ such that $M^*$ satisfies: $M^*=8\,\pi$ if $0<\tau\le\tau^*$ and $M^*>8\,\pi$ if $\tau>\tau^{**}$.
\label{the:existence2}
\end{theorem}
When $M^*>8\,\pi$, there are at least two positive solutions for any $M\in(8\,\pi,M^*)$. When $M^*=8\,\pi$, it is still an open question to decide if there is a positive solution $(\phi,S)\in C^2[0,\infty)\times C^1[0,\infty)$ such that $M=M^*$ or to prove a uniqueness result for any $M\in(0,8\,\pi)$.
\begin{remark}
\label{mono} \
The estimate $\tau^*>1/2$ will be given in Proposition~\ref{lem:estM1}, as well as refined estimates on $M(a,\tau)$. Theoretical results show that $\tau^*\in(0.5,16.1109\ldots)$, see Th.~\ref{the:existence}, \eqref{f}--\eqref{g} and Sec.~\ref{Sec:Estimate}, while numerical computations suggests that $\tau^*\in(0.62,0.64)$, see Fig. 2 (right). Moreover, it is an interesting open question to decide whether $\tau^*=\tau^{**}$, as again the numerical results suggest, or not. Exact multiplicities of solutions for $M>8\,\pi$ are not known in detail either. Let us observe that for $\tau>\tau^*$ the function $M(a,\tau)$ depends on $a$ in a nonmonotone manner. This is a significant difference with the monotone dependence of self-similar solutions of the parabolic-elliptic Keller--Segel system (see \cite[Sec. 4]{BKLN}).
\end{remark}

\section{Qualitative properties of $\phi$ and $S$}
\label{Sec:qualitative}

In the present section we will derive all {\em a priori} estimates on $\phi$ and $S$ which are necessary to prove Theorems \ref{the:existence}, \ref{the:BU} and \ref{the:existence2}. Some of them are new while other were already known. In any case, we shall give a unified and simplified proof of all of them in terms of cumulated densities.

\subsection{Preliminary estimates}

Let $(u,v)\in(C^2_0(\R^2))^2$ be a positive solution of \eqref{KSuSS}--\eqref{KSvSS} with $u\in L^1(\R^2)$.
The corresponding $(\phi,S)$ satisfies \eqref{phi2}--\eqref{S}, \eqref{IC} with $a=u(0)/2$. Moreover, for any $y>0$, it immediately holds true that: $\phi$ is a positive, strictly increasing and concave function on $(0,\infty)$ while $0<S(y)<\phi(y)$ for any $y>0$ since $S'<\phi'$ on $(0,\infty)$ by \eqref{S}. More precisely, an integration by parts in \eqref{S1} gives
\begin{equation}
S(y)=\phi(y)-\frac\tau4\,{\rm e}^{-\tau\,y/4}\int_0^y{\rm e}^{\tau\,z/4}\,\phi(z)\dz\,.
\label{S2}
\end{equation}
On the other hand, in \eqref{S2}, the increasing monotonicity property of $\phi$ gives us
\begin{equation}
S(y)\ge\phi(y)-\frac\tau4\,{\rm e}^{-\tau\,y/4}\,\phi(y)\int_0^y{\rm e}^{\tau\,z/4}\dz= {\rm e}^{-\tau\,y/4}\,\phi(y)\,,
\label{estS1}
\end{equation}
while the decreasing monotonicity property of $\phi'$ in \eqref{S1} leads to
\begin{equation}
S(y)\ge {\rm e}^{-\tau\,y/4}\,\phi'(y)\int_0^y{\rm e}^{\tau\,z/4}\dz=\frac4\tau\,\phi'(y)\left(1-{\rm e}^{-\tau\,y/4}\right)\;
\label{estS2}
\end{equation}
for each $y\ge 0$. From \eqref{estS1} and \eqref{estS2}, we get
\[
\frac\tau2\,S(y)\ge\left(\phi(y)-\phi(y)\,{\rm e}^{-\tau\,y/4}\right)'.
\]
Since $\frac\tau2\,S=2\,\phi'-2\,S'$, the last inequality gives
\begin{equation*}
S(y)\le\frac12\,\phi(y)\left(1+ {\rm e}^{-\tau\,y/4}\right)
\label{estS3}
\end{equation*}
for each $y\ge 0$, which is a better estimate than $S<\phi$ but still not yet satisfactory for large $y$.
\medskip

Let us now estimate $\phi$. Looking closer at system \eqref{phi2}--\eqref{S}, one observes that the quantity ${\rm e}^{\,y/4}\,\phi'(y)$ is positive and decreasing. Hence
\begin{equation}
l(a,\tau):=\lim_{z\to\infty}{\rm e}^{z/4}\;\phi'(z)\le {\rm e}^{\,y/4}\;\phi'(y)\le \phi'(0)=a\,,
\label{estphi'}
\end{equation}
for any $y\ge 0$. Notice that $l(a,\tau)=\sigma/2$, which proves that $\lim_{s\to-\infty}\sigma(s)=0$ with the notations of  Remark~\ref{rk:sigma}. Integrating once more the above inequalities on $[0,y]$ we have
\begin{equation}
4\,l(a,\tau)\left(1-{\rm e}^{-y/4}\right)\le\phi(y)\le4\,a\left(1-{\rm e}^{-y/4}\right)\,.
\label{estphi}
\end{equation}
In particular, for each $\tau>0$, $M(a,\tau)$ is finite,
\begin{equation}
l(a,\tau)\le\frac{M(a,\tau)}{8\,\pi}\le a\,,
\label{estM}
\end{equation}
and we see that, whatever $\tau$ is, the shooting parameter $a$ has to be large enough ($a>1$) in order to obtain a self-similar solution $u$ with mass $M>8\,\pi$.

We can improve estimate \eqref{estphi} as follows. Since $\lim_{y\to\infty}\phi'(y)=0$, integrating the inequality $\phi''+\frac14\,\phi'<0$ on $[y,\infty)$, we get
\[
\phi'(y)+\frac14\,\phi(y)\ge\frac {M(a,\tau)}{8\,\pi}\,,
\]
and therefore, by integrating once more on $[0,y]$,
\[
\phi(y)\ge\frac {M(a,\tau)}{2\,\pi}\(1-{\rm e}^{-y/4}\)\,.
\]
In conclusion, using the previous estimate for $\phi$, we obtain for each $y\ge 0$
\begin{equation}
\frac {M(a,\tau)}{2\,\pi}\left(1-{\rm e}^{-y/4}\right) \le\phi(y)\le\min\left\{4\,a\left(1-{\rm e}^{-y/4}\right),\frac {M(a,\tau)}{2\,\pi}\right\}\,,
\label{estphi2}
\end{equation}
where equality in the minimum is achieved for $\overline y=-\,4\log\big(1-\frac M{8\,\pi\,a}\big)\in(0,\infty]$. In particular, equalities hold in \eqref{estphi2}, {\em i.e.} $\phi(y)=\frac M{2\,\pi}\big(1-{\rm e}^{-y/4}\big)$, if and only if $M=8\,\pi\,a$, in which case $\overline y=\infty$. But since $\phi(y)=\frac M{2\,\pi}\big(1-{\rm e}^{-y/4}\big)$ is not a solution of \eqref{phi2}--\eqref{S}, estimate~\eqref{estphi2} holds true with strict inequalities as well as $M<8\,\pi\,a$.
\medskip

Coming back to the function $S$, using estimate \eqref{estphi'} and identity \eqref{S1}, we have
\[
S(y)\le a\;{\rm e}^{-\tau\,y/4}\int_0^y{\rm e}^{(\tau-1)\,z/4}\dz\,,
\]
{\em i.e.}, for each $y\ge 0$ and $\tau>0$,
\begin{equation}
S(y)\le a\;y\;h(y;\tau)
\label{estStau}
\end{equation}
where
\begin{equation}
h(y;\tau)=\left\{
\begin{array}{rcll}
&&{\rm e}^{-y/4}\quad\mbox{if}\;\tau=1\,,\\
&&\frac{4}{y\,(\tau-1)}\({\rm e}^{-y/4}-{\rm e}^{-\tau\,y/4}\)\quad\mbox{if}\;\tau\neq1\,.
\end{array}
\right.
\label{eq:h}
\end{equation}
As a consequence, it holds true that
\[
\lim_{y\to\infty}S(y)=\lim_{y\to\infty}\frac{S(y)}y=0\,,
\]
$S(y)/y$ is integrable near $y=0$ and, using \eqref{estS2}, $S'(0)=a$.

The above asymptotic behavior of $S$ at infinity, together with the initial condition $S(0)=0$, allow us to integrate equation \eqref{S} on $[0,\infty)$ to obtain
\begin{equation}
\frac {M(a,\tau)}{2\,\pi}=\phi(\infty)=\frac\tau4\int_0^{\infty}S(y)\dy\,.
\label{M}
\end{equation}
Therefore, any appropriate bound for $S$ would give a bound for the total mass $M$. However, let us observe that if we plug estimates \eqref{estStau} into \eqref{M}, we found again the upper bound in \eqref{estM}. Finally, thanks to the integrability of $S(y)/y$ near $y=0$, equation~\eqref{phi2} written as
\[
\phi''+\phi'\left(y/4+\frac12\int_0^y\frac{S(z)}z\dz\right)'=0
\]
and integrated on $[0,y]$ gives the relation
\begin{equation}
\phi'(y)=a\;{\rm e}^{-y/4}\;\exp\left(-\frac12\int_0^y\frac{S(z)}z\dz\right)\,.
\label{formule}
\end{equation}
%

\subsection{Further estimates}

First, let us improve on the lower bound in~\eqref{estphi2} for $\phi$. As far as we know, all estimates of this section are new. Using the fact that $S<\phi$ in~\eqref{phi2}, for $y>0$ we have
\[
\phi''+\frac14\,\phi'+\frac1{2\,y}\,\phi'\phi>0\,.
\]
After a multiplication by $y$, an integration on $[0,y]$ leads to
\[
y\,\phi'-\phi+\frac y4\,\phi+\frac14\,\phi^2>\frac14\int_0^y\phi(z)\dz\,.
\]
Dividing by $\phi^2\;{\rm e}^{\,y/4}$ we obtain the differential inequality
\[
\left(-\frac y\phi\,{\rm e}^{-y/4}\right)'+\frac14\,{\rm e}^{-y/4}>\frac14\,\frac1{\phi^2}\,{\rm e}^{-y/4}\int_0^y\phi(z)\dz\,.
\]
Finally, dropping the positive term on the right hand side, and integrating once again on $[0,y]$ gives us a lower bound for any $\tau>0$ and $a>0$, namely,
\[
\phi(y)\ge\frac y{\(1+\frac1a\)\,{\rm e}^{\,y/4}-1}
\]
for each $y\ge 0$. 
This is, of course, a better estimate than \eqref{estphi2} but only for $y$ near the origin since the inequality $S(y)<\phi(y)$ is a good approximation for $y$ near the origin but not for large~$y$. However, we can now replace \eqref{estphi2} with
\begin{equation}
\max\left\{{\textstyle\tfrac {M(a,\tau)}{2\,\pi}\left(1-{\rm e}^{-y/4}\right),\tfrac y{\(1+\frac1a\)\,{\rm e}^{\,y/4}-1}}\right\} \le\phi(y)\le\min\left\{{\textstyle 4\,a\left(1-{\rm e}^{-y/4}\right)\,,\tfrac {M(a,\tau)}{2\,\pi}}\right\}.
\label{estphi3}
\end{equation}
The maximum on the left hand side of \eqref{estphi3} is achieved by both terms at some $\tilde y>0$ and
\[
\max\left\{\textstyle{\frac{M(a,\tau)}{2\,\pi}\left(1-{\rm e}^{-y/4}\right)\,,\frac y{\(1+\frac1a\){\rm e}^{\,y/4}-1}}\right\}=\frac y{\(1+\frac1a\)\,{\rm e}^{\,y/4}-1}
\] 
for each $y\in[0,\tilde y]$. 
Moreover, for any $y\ge y^*$, we have
\[
\frac{M(a,\tau)}{2\,\pi}>\phi(y)\ge\phi(y^*)\ge \frac {y^*}{\(1+\frac1a\)\,{\rm e}^{\,y^*/4}-1} =\frac 4{1+\frac1a}\;\to 4^-\quad \mbox{as}\; a\to\infty\,,
\]
if $y^*$ is the point where the maximum of $y\mapsto\frac y{\(1+\frac1a\)\,{\rm e}^{\,y/4}-1}$ is achieved.

Next, let us apply estimates \eqref{estStau}--\eqref{eq:h} to \eqref{formule}. For $\tau\neq1$, we have
\begin{equation*}
\begin{array}{rcl}
\displaystyle\int_0^yh(z;\tau)\dz&&\displaystyle=\frac4{\tau-1}\int_0^y\frac1z\int_\tau^1\frac{\rm d}{{\rm d}t}({\rm e}^{-\frac t4z})\dt\dz
=\frac1{\tau-1}\int_0^y\int_1^\tau {\rm e}^{-\frac t4z}\dt\dz
\\ \\
&&\displaystyle=\frac4{\tau-1}\int_1^\tau \frac1t\(1-{\rm e}^{-t\,y/4}\)\dt
=\frac4{\tau-1}\log\tau-\frac4{\tau-1}\int_1^\tau \frac1t\,{\rm e}^{-t\,y/4}\dt\,,
\end{array}
\end{equation*}
\[
\int_0^y\frac{S(z)}z\dz\le 4\,a\,I(\tau)
\]
with $I(\tau)=\frac{\log\tau}{\tau-1}$ and
\begin{equation}
\phi'(y)\ge a\,{\rm e}^{-y/4}{\rm e}^{-2\,a\,I(\tau)}\,.
\label{estphi'4}
\end{equation}
Integrating \eqref{estphi'4} on $[0,\infty)$, we get the same estimate as in \eqref{est:Mgrand} giving arbitrarily large mass $M$ for $\tau$ large enough, {\em i.e.}
\begin{equation}
\frac {M(a,\tau)}{2\,\pi}\ge 4\,a\,{\rm e}^{-2\,a\,I(\tau)}\,.
\label{estM4}
\end{equation}
For $\tau=1$, since $h(y;1)={\rm e}^{-y/4}$, one obtains, for all $a>0$,
\begin{equation*}
\frac {M(a,1)}{2\,\pi}\ge2\(1-{\rm e}^{-2\,a}\)\,.
\label{estM5}
\end{equation*}
%
\begin{remark}\label{compa} \
The lower bound \eqref{estM4} is compatible with the upper bounds for $M$ known from \cite{B06}, {\em i.e.}
\begin{multline*}
\frac M{2\,\pi}\le 4\;\mbox{if}\;\tau\in\left(0,1/2\right]\,,\quad
\frac M{2\,\pi}\le \frac23\,\pi^2\;\mbox{if}\;\tau\in\left(1/2,1\right]\\\mbox{and}\quad\frac M{2\,\pi}\le \min\left\{{\textstyle\frac23\,\pi^2\,\tau,4\,(\tau+1)}\right\}\;\mbox{if}\;\tau>1\,.
\end{multline*}
\end{remark}

Finally, following \cite{NSY02}, define the new function
\[
W(y):=\int_0^y\phi'(z)\,{\rm e}^{\tau\,z/4}\dz={\rm e}^{\tau\,y/4}\,S(y)\,,
\]
where the second equality follows from \eqref{S1}. After a multiplication of \eqref{phi3} by ${\rm e}^{\tau\,y/4}$, it is easy to see that $W$ satisfies the initial value
problem
\begin{eqnarray*}
&&W''+\frac{1-\tau}4\,W'+\frac1{4\,y}\(W^2\)'{\rm e}^{-\tau\,y/4}=0\,,\\
&&W(0)=0\,,\quad W'(0)=a\,.
\end{eqnarray*}
Next, a multiplication by $y$ and an integration on $[0,y]$ gives us
\[
y\,W'-W+\frac{1-\tau}4\int_0^yz\,W'(z)\dz+\frac14\,{\rm e}^{-\tau\,y/4}\,W^2+\frac\tau{16}\int_0^y{\rm e}^{-\tau\,z/4}\,W^2(z)\dz=0\,.
\]
Dividing by $W^2$ the equation becomes
\begin{equation}
\left(-\frac yW\right)'+\frac{1-\tau}4\frac1{W^2}\int_0^yz\,W'(z)\dz+\frac14\,{\rm e}^{-\tau\,y/4}+\frac\tau{16}\frac1{W^2}\int_0^y{\rm e}^{-\tau\,z/4}\,W^2(z)\dz=0\,.
\label{eqnW1}
\end{equation}
This last identity is a useful reformulation of the problem for $0<\tau\le1$, since in this case the two integral terms in the equation are positive. Then, eliminating both of them and integrating on $[0,y]$, we get for each $y\ge 0$
\[
\frac y{W(y)}\ge \frac1a+\frac1\tau\left(1-{\rm e}^{-\tau\,y/4}\right),
\]
{\em i.e.}
\begin{equation}
S(y)\le\frac{\tau\,a\,y}{\tau\,{\rm e}^{\tau\,y/4}+a\({\rm e}^{\tau\,y/4}-1\)}\,.
\label{estStau<=1}
\end{equation}
For $\tau>1$ it is more convenient to integrate by parts the first integral term in \eqref{eqnW1} to obtain
\begin{multline*}
\left(-\frac yW\right)'+\frac{1-\tau}4\frac y{W}-\frac{1-\tau}4\frac1{W^2}\int_0^yW(z)\dz+\frac14\,{\rm e}^{-\tau\,y/4}\\
+\frac\tau{16}\frac1{W^2}\int_0^y{\rm e}^{-\tau\,z/4}\,W^2(z)\dz=0\,.
\label{eqnW2}
\end{multline*}
Again, eliminating the two positive integral terms and multiplying by ${\rm e}^{(\tau-1)\,y/4}$, we obtain
\begin{equation*}
\left({\rm e}^{(\tau-1)\,y/4}\frac yW\right)'\ge\frac14\,{\rm e}^{-y/4}.
\end{equation*}
After an integration on $[0,y]$, this gives
\begin{equation*}
{\rm e}^{(\tau-1)\,y/4}\frac yW\ge\frac1a+1-{\rm e}^{-y/4}\,,
\end{equation*}
{\em i.e.}
\begin{equation}
S(y)\le\frac{a\,y}{{\rm e}^{\,y/4}+a\({\rm e}^{\,y/4}-1\)}\,.\label{estStau>1}
\end{equation}
Summarizing, estimates \eqref{estStau<=1} and \eqref{estStau>1} read
\begin{equation}
S(y)\le a\;y\;g(y;a,\tau)\,,
\label{estStau2}
\end{equation}
for all $y\ge 0$, $a>0$, $\tau>0$,
where
\begin{equation}
g(y;a,\tau):=\left\{
\begin{array}{rcll}
&&\displaystyle\frac\tau{(\tau+a)\,{\rm e}^{\tau\,y/4}-a}\quad&\mbox{if}\;0<\tau\le1\,,\\ \\
&&\displaystyle\frac1{(1+a)\,{\rm e}^{\,y/4}-a}\quad&\mbox{if}\;\tau>1\,.
\end{array}
\right.
\label{eq:g}
\end{equation}
As an important consequence of \eqref{estStau2}--\eqref{eq:g}, for any $\tau>0$, $S$ is bounded uniformly with respect to $a>0$:
\begin{equation}
S(y)\le\frac{\min\{\tau,1\}\,y}{{\rm e}^{\min\{\tau,1\}\,y/4}-1}
\label{eq:SUnifBded}
\end{equation}%
for each $y>0$. 
Such an estimate does not follow from \eqref{estStau}--\eqref{eq:h}.

Estimate \eqref{estStau2} is better than estimate \eqref{estStau} for $\tau=1$. For $\tau\neq1$, this depends on the values of $\tau$ and $a$. Therefore, it is interesting to reproduce the computations giving~\eqref{estM4} by using the function $g$ instead of $h$. For $\tau\ge1$ and each $y\ge 0$, we obtain
\[
\int_0^yg(z;a,\tau)\dz=\frac4a\,\log\left[(1+a)\,{\rm e}^{\,y/4}-a\right]-\frac ya\,,
\]
and from equation \eqref{formule}
\begin{equation}
\phi'(y)\ge a\,\frac{{\rm e}^{\,y/4}}{\left[(1+a)\,{\rm e}^{\,y/4}-a\right]^2}\,.
\label{estphi'5}
\end{equation}
This gives, for $a>0$ and $\tau\ge 1$,
\begin{equation}
\frac {M(a,\tau)}{2\,\pi}\ge\frac{4\,a}{a+1}\,.
\label{estM6}
\end{equation}
Such a lower bound is definitely worse than \eqref{estM4} for large values of $\tau$ or, to be precise, as soon as $I(\tau)\le\log(a+1)/(2\,a)$. On the other hand, for $\tau<1$, we have
\[
\int_0^yg(z;a,\tau)\dz=\frac4a\log\left[\left(\frac a\tau+1\right){\rm e}^{\tau\,y/4}-\frac a\tau\right]-\frac\tau a\,y\,,
\]
and again from equation \eqref{formule},
\begin{equation*}
\phi'(y)\ge a\,{\rm e}^{-y/4}\,{\rm e}^{\tau\frac y2}\frac1{\left[\(\frac a\tau+1\){\rm e}^{\tau\,y/4}-\frac a\tau\right]^2}\,.
\label{estphi'6}
\end{equation*}
Finally, it holds true that, for $a>0$ and $\tau<1$,
\begin{equation}
\frac {M(a,\tau)}{2\,\pi}\ge\frac{4\,a\,\tau}{a+\tau}\,.
\label{estM7}
\end{equation}
To conclude, integrating \eqref{estphi'5} on $[0,y]$ and using estimate \eqref{estStau2} gives us, for any $\tau\ge1$,
\[
S(y)\le\frac y{(1+\frac1a)\,{\rm e}^{\,y/4}-1}
\le\frac {4\,({\rm e}^{\,y/4}-1)}{(1+\frac1a)\,{\rm e}^{\,y/4}-1}\le\phi(y)
\]
for each $y\ge 0$, 
which is a good approximation of $S$ and $\phi$ near the origin since it takes into account the condition $S'(0)=\phi'(0)=a$. Moreover, \eqref{estphi3} is improved and replaced with
\[
\max\left\{{\textstyle\frac {M(a,\tau)}{2\,\pi}\left(1-{\rm e}^{-y/4}\right),\frac {4\,({\rm e}^{\,y/4}-1)}{\(1+\frac1a\){\rm e}^{\,y/4}-1}}\right\}\le\phi(y)\le\min\left\{{\textstyle4\,a\left(1-{\rm e}^{-y/4}\right),\frac{M(a,\tau)}{2\,\pi}}\right\}
\]
for any $\tau\ge1$ and $y\ge 0$.

\begin{remark}\label{rk:sigma(s)} \
As an additional consequence of the above estimates, we observe that
\[
a\int_0^yg(z;a,\tau)\dz=4\,\log\left(1+\tfrac a{\min\{1,\tau\}}\,-\tfrac a{\min\{1,\tau\}}e^{-\min\{1,\tau\}\,y/4}\,\right)
\]
converges as $y\to\infty$ , so that
\[
\exp\left[-\frac a2\int_0^\infty g(z;a,\tau)\dz\right]=\(1+\tfrac a{\min\{1,\tau\}}\)^{-2}\,.
\]
According to \eqref{estphi'}, \eqref{formule} and \eqref{estStau2}, we find the estimate
\[
\tfrac\sigma2=l(a,\tau)=\lim_{y\to+\infty}a\,\exp\left(-\frac12\int_0^y\frac{S(z)}z\dz\right)
\ge a\(1+\tfrac a{\min\{1,\tau\}}\)^{-2}\,,
\]
which, taking into account the change of parametrization $s=\log(2\,a)$, refines the estimate $\lim_{s\to+\infty}\sigma(s)=0$ found in \cite{NSY02} and our estimate \eqref{est:sigma}. Here we use the notations of Remark~\ref{rk:sigma}.
\end{remark}

\subsection{New upper bounds}

Using the previous estimates on $S$ and an argument in \cite{B06}, we can improve on the upper bound in \eqref{estM}. Let
\begin{equation*}
j(\tau):=\left\{
\begin{array}{ccl}
\infty\quad&\mbox{if}&\;0<\tau\le\frac12\,,\\ \\
\tau\frac{{\rm e}^{1-\frac1{2\,\tau}}}{2\,\tau-{\rm e}^{1-\frac1{2\,\tau}}}\quad&\mbox{if}&\;\frac12<\tau\le1\,,\\ \\
\frac{{\rm e}^{1-\frac1{2\,\tau}}}{2\,\tau-{\rm e}^{1-\frac1{2\,\tau}}}\quad&\mbox{if}&\;\tau>1\,.
\end{array}
\right.
\end{equation*}
%
\begin{proposition}
For any $\tau>0$, if $a\le\max\{j(\tau),1\}$, then $M(a,\tau)\le8\,\pi\min\{1,a\}$.
\label{lem:estM1}
\end{proposition}
The above estimate gives us a nonoptimal set of parameters $(a,\tau)$ that guarantees $M(a,\tau)\le8\,\pi$. It is interesting to notice that $\lim_{\tau\to(1/2)^+}j(\tau)=\infty$.
\begin{proof} Let $M=M(a,\tau)$. From the identity
\[
\left(\frac M{2\,\pi}\right)^2-4\left(\frac M{2\,\pi}\right)=\int_0^{\infty}\left(2\,\phi(y)\,\phi'(y)+4\,y\,\phi''(y)\right) \dy\,,
\]
and $4\,y\,\phi''=-y\,\phi'-2\,\phi' S$ which follows from \eqref{phi2}, we have, after an integration by parts and using \eqref{S},
\begin{multline*}
\left(\frac M{2\,\pi}\right)^2-4\left(\frac M{2\,\pi}\right)=\int_0^{\infty}\phi'(2\,\phi-2\,S-y)\dy
=\int_0^{\infty}\left(\phi-\frac M{2\,\pi}\right)'(2\,\phi-2\,S-y)\dy\\
=-\int_0^{\infty}\left(\phi-\frac M{2\,\pi}\right)(2\,\phi'-2\,S'-1)\dy
=-\int_0^{\infty}\left(\phi-\frac M{2\,\pi}\right)\left(\frac\tau2\,S-1\right)\dy\,.
\end{multline*}
Hence we have $\frac M{2\,\pi}\le4$ if
\begin{equation}
\frac\tau2\,S(y)\le1
\label{SC}
\end{equation}
for each $y>0$. 
{}From \eqref{eq:SUnifBded} it follows that $S(y)<4$ for all $y\ge0$, for any $\tau>0$ and $a>0$: the above sufficient condition \eqref{SC} is satisfied whenever $\tau\le1/2$. For $\tau>1/2$ we have to use one of the previous upper bounds for $S$.

$(a)$ Using \eqref{estStau}, we have
\[
1-\frac\tau2\,S(y)\ge1-2\,a\,\frac\tau{\tau-1}\left({\rm e}^{-y/4}-{\rm e}^{-\tau\,y/4}\right)
\]
for any $\tau\neq1$ and each $y\ge 0$, and condition \eqref{SC} is satisfied if
\[
a\le\min_{y>0}\frac12\frac{\tau-1}\tau\frac1{{\rm e}^{-y/4}-{\rm e}^{-\tau\,y/4}}=\frac12\,\tau^{\frac{1}{\tau-1}}\,.
\]
For $\tau=1$, using \eqref{estStau} as before (or by continuity of the previous argument as $\tau\to1$), we similarly obtain
\[
a\le\min_{y>0}\frac{2\,{\rm e}^{\,y/4}}y=\frac{\rm e}2\,.
\]

$(b)$ Using \eqref{estStau2}, we have for $\tau>1$ and each $y\ge 0$
\[
\frac\tau2\;S(y)-1\le\frac12\;\frac{\tau\,a\,y}{(1+a)\,{\rm e}^{\,y/4}-a}-1\,,
\]
Then condition \eqref{SC} is satisfied if
\[
a\le\min_{0<y<\bar{y}}\;\frac{2\,{\rm e}^{\,y/4}}{\tau\,y-2\,({\rm e}^{\,y/4}-1)}
=\frac{{\rm e}^{1-\frac1{2\,\tau}}}{2\,\tau-{\rm e}^{1-\frac1{2\,\tau}}}\,,
\]
where we take into account that $\tau\,y-2\,({\rm e}^{\,y/4}-1)<0$ for $y>\bar{y}$, $\bar{y}$ being the unique solution of the equation $\frac\tau2\,y+1={\rm e}^{\,y/4}$. Similarly, for $\frac12<\tau\le1$ and each $y\ge 0$, we get
\[
\frac\tau2\;S(y)-1\le\frac\tau2\;\frac{\tau\,a\,y}{(\tau+a)\,{\rm e}^{\frac \tau4\,y}-a}-1\,.
\]
Then condition \eqref{SC} is satisfied if
\[
a\le\min_{0<\tau\,y<\bar{y}}\;\frac{2\,\tau\,{\rm e}^{\frac \tau4\,y}}{\tau^2\,y-2\({\rm e}^{\frac \tau4\,y}-1\)}
=\tau\frac{{\rm e}^{1-\frac1{2\,\tau}}}{2\,\tau-{\rm e}^{1-\frac1{2\,\tau}}}\,.
\]
Comparing the results obtained in $(a)$ and $(b)$, the proof of Proposition~\ref{lem:estM1} is completed.
\qed
\end{proof}

\section{Proofs}
\label{Sec:Proofs}

This section is devoted to the proof of Theorems \ref{the:existence}, \ref{the:BU} and \ref{the:existence2}. As a byproduct of these results, we obtain Theorem \ref{Thm:Simplified}.

\subsection{Proof of Theorem \ref{the:existence}}

Given any fixed $(a,\tau)\in\R^2_+\,$, the local existence issue of the (singular) system \eqref{phi2}--\eqref{S} with initial conditions \eqref{IC} can be solved using a fixed point argument applied to the operator
\[
{\mathcal T}[\Phi](y)=a\,{\rm e}^{-y/4}-\frac12\,{\rm e}^{-y/4}\int_0^y\frac1z\,{\rm e}^{(1-\tau)\,z/4}\,\Phi(z)\int_0^z{\rm e}^{\tau\,\xi/4}\,\Phi(\xi)\,\dxi\,,
\]
defined on the complete metric space $X_a:=\{\Phi\!\in\!C[0,y_a]\,:\,\Phi(0)=a,\,$\hbox{$0\le\Phi(y)\le a$}, $\ 0\le y\le y_a\}$ endowed with the usual supremum norm. Indeed, an appropriate choice of $y_a$ gives that $\mathcal T$ maps $X_a$ into $X_a$ and that $\mathcal T$ is a contraction. If ${\mathcal T}[\Phi]=\Phi$, it is then enough to define $\phi(y):=\int_0^y\Phi(z)\,\dz$ and $S(y):= {\rm e}^{-\tau\,y/4}\,\int_0^y{\rm e}^{\tau\,z/4}\,\Phi(z)\dz$ in order that $(\phi,S)$ is a solution of \eqref{phi2}--\eqref{S}, \eqref{IC} with $\phi\in C^1[0,y_a]\cap C^2(0,y_a]$ and $S\in C^1[0,y_a]$. The continuation of the local solution to a global one is standard since system \eqref{phi2}--\eqref{S} is no more singular away from the origin and solutions are locally bounded on $\R_+$ by the estimates of Section~\ref{Sec:qualitative}.

The fact that $\phi\in C^2[0,\infty)$ follows from \eqref{formule} and $\lim_{y\to0_+}S(y)/y=a=S'(0)$. Estimates \eqref{f} have been proved in Section~\ref{Sec:qualitative}.

Finally, uniqueness of global solutions of \eqref{phi2}--\eqref{S}, \eqref{IC} is a consequence of the contraction property of $\mathcal T$ and 
the Cauchy--Lipschitz theorem.

Concerning the continuity of the map $a\in\R_+\mapsto (\phi,S)$, let us denote by $(\phi_i,S_i)$ the solution associated to the shooting parameter $a_i$, $i=1,2$. Following \cite{MMY99} we have
\begin{equation}
|\log\phi_1'(y)-\log\phi_2'(y)|\le |\log a_1-\log a_2\,|+\frac12\int_0^y\frac1z\,|S_1(z)-S_2(z)|\dz
\label{cont1}
\end{equation}
and
\begin{multline}
|S_1(y)-S_2(y)|\le {\rm e}^{-\tau\,y/4}\int_0^y{\rm e}^{(\tau-1)\,z/4}\left|\,{\rm e}^{z/4}\phi'_1(z)-{\rm e}^{z/4}\phi'_2(z)\,\right|\dz\\
\le {\rm e}^{\max\{\log a_1,\log a_2\}}\,{\rm e}^{-\tau\,y/4}\int_0^y{\rm e}^{(\tau-1)\,z/4}\,|\,\log\phi'_1(z)-\log\phi'_2(z)\,|\dz\,,
\label{cont2}
\end{multline}
where the decreasing monotonicity property of the function ${\rm e}^{\,y/4}\,\phi'(y)$ has been used in the last inequality. Plugging \eqref{cont2} into \eqref{cont1} and denoting $C={\rm e}^{\max\{\log a_1,\log a_2\}}$, we obtain
\begin{multline}
|\log\phi_1'(y)-\log\phi_2'(y)|\\
\le |\log a_1-\log a_2\,|+\frac C2\int_0^y\frac1z\,{\rm e}^{-\tau\,z/4}\int_0^z{\rm e}^{(\tau-1)\,\zeta/4}|\log\phi'_1(\zeta)-\log\phi'_2(\zeta)|\dzeta \dz\\
\le |\log a_1-\log a_2\,|+\frac C2\int_0^y|\log\phi'_1(\zeta)-\log\phi'_2(\zeta)\,|\,f(\zeta)\dzeta\,,
\label{cont3}
\end{multline}
where $f(\zeta)={\rm e}^{(\tau-1)\,\zeta/4}\int_\zeta^{\infty}\frac1z\, {\rm e}^{-\tau\,z/4}\dz$. Next, $f\in L^1(0,\infty)$ with $\int_0^{\infty}f(\zeta)\dzeta=4\,\frac{\log\tau}{\tau-1}=4\,I(\tau)$. Therefore, the Gronwall lemma applied to \eqref{cont3} gives us
\begin{equation}
|\log\phi_1'(y)-\log\phi_2'(y)|\le |\log a_1-\log a_2\,|\,{\rm e}^{\frac C2\!\int_0^yf(\zeta)\dzeta}\le |\log a_1-\log a_2\,|\,{\rm e}^{2\,C\,I(\tau)}.
\label{cont4}
\end{equation}
Estimate \eqref{cont4} implies the continuity of the map $a\mapsto\phi'$. The continuity of the maps $a\mapsto S$ and $a\mapsto\phi$ follows by \eqref{cont2}--\eqref{cont4} and by the identity $\phi(y)=S(y)+\frac\tau4\int_0^yS(z)\dz$ respectively. Finally, the continuity of $a\mapsto M$ follows by $\frac M{2\,\pi}=\frac\tau4\int_0^{\infty}S(y)\dy$; see Section~\ref{Sec:qualitative} for more details.
\qed

\subsection{Proof of Theorem \ref{the:BU}}

The existence of a sequence of positive self-similar solutions $(u_k,v_k)$ corresponding to a positive sequence $\{a_k\}$ is an immediate consequence of the existence of a positive solution $(\phi_k,S_k)$ of \eqref{phi2}--\eqref{S}, \eqref{IC} by Theorem \ref{the:existence}. Indeed, it is sufficient to define
\[
u_k(r)=2\,\phi_k'(r^2)\quad\mbox{and}\quad v_k(r)=\frac12\int_ {r^2}^\infty\frac{S_k(z)}z\dz\,,
\]
as follows from \eqref{Derphi} and \eqref{eq:v(0)}. Moreover, $u_k\in C^1[0,\infty)$ and $v_k\in C^2[0,\infty)$. Whenever $a_k\to\infty$, the limit $\|v_k\|_{L^\infty(\R^2)}\to\infty$ follows from $\|v_k\|_{L^\infty(\R^2)}=v_k(0)$ and \eqref{est3:v(0)}.

Next, let us define $M_k:=\|u_k\|_{L^1(\R^2)}$. From the estimates of Section~\ref{Sec:qualitative}, the sequence $\{M_k\}$ is bounded from above (by a constant depending on $\tau$), and there exist two subsequences, still denoted $M_k$ and $u_k$, such that $M_k\to\alpha$ and $u_k\rightharpoonup\,\alpha\,\delta_0$. The delta measure is centered at $\xi=0$ since $u_k(0)=2\,a_k$. Actually $\alpha=8\,\pi$ for any $\tau>0$, as an immediate consequence of the identity obtained in the proof of Proposition \ref{lem:estM1}
\[
\left(\frac{M_k}{2\,\pi}\right)^2-4\left(\frac{M_k}{2\,\pi}\right)=\int_0^{\infty}\phi_k'(2\,\phi_k-2\,S_k-y)\dy\,.
\]
Hence we have
\begin{equation}
\left(\frac {M_k}{2\,\pi}\right)^2-4\left(\frac {M_k}{2\,\pi}\right)=\frac1\pi\,\int_{\R^2} u_k(\xi)\left(\phi_k(|\xi|^2)-S_k(|\xi|^2)-\tfrac12|\xi|^2\right)\dxi\,.
\label{eq:DiracAsymp}
\end{equation}
Letting $k\to\infty$ and observing that:
\begin{description}
\item[(i)] $(\phi_k-S_k)'=\tau\,S_k/4$ is bounded by \eqref{eq:SUnifBded}, uniformly with respect to $a_k\to\infty$,
\item[(ii)] $u_k(r)=2\,\phi_k'(r^2)$ is uniformly decaying (with respect to $a_k\to\infty$) for large values of $r$ and $\lim_{r\to\infty}\sup_k\int_{|\xi|>r} u_k(\xi)\,|\xi|^2\dxi=0$, as a consequence of \eqref{estS2} and \eqref{eq:SUnifBded},
\end{description}
we obtain that the right hand side in \eqref{eq:DiracAsymp} converges to $0$. On the other hand, $\alpha$ is necessarily positive by \eqref{estM6} and \eqref{estM7}, which proves that $\alpha=8\,\pi$.
\qed

\begin{remark}
\label{rk:BU} \
Let us observe that the identity
\[
4\,M+2\int_{\R^2} u(\xi)\,\nabla v(\xi)\cdot\xi\,\dxi-\int_{\R^2}|\xi|^2\,u(\xi)\dxi=0
\]
follows from equation \eqref{KSuSS} multiplied by $|\xi|^2$ and from the integrability of $u$ given by \eqref{u} and \eqref{est:v}. Mimicking a standard computation for the parabolic-elliptic Keller--Segel system by writing $v=-\frac1{2\,\pi}\log(\cdot)*u+\tilde v$, the above identity reads
\[
4\,M-\frac{M^2}{2\,\pi}+2\int_{\R^2} u(\xi)\,\nabla \tilde v(\xi)\cdot\xi\,\dxi-\int_{\R^2}|\xi|^2 u(\xi)\dxi=0\,.
\]
See for instance \cite{BDP,MR2103197} for more details. Therefore, we have found that $\nabla\tilde v(\xi)\cdot\xi=\phi(|\xi|^2)-S(|\xi|^2)\ge 0$. This is consistent with the fact that, from equation \eqref{KSvr}, one easily finds that $\phi(r^2)-S(r^2)=-\frac\tau2\int_0^rs^2\,v'(s)\ds$.
\end{remark}

\subsection{Proof of Theorem \ref{the:existence2}}

For any fixed $\tau$, let us define $M^*(\tau)=\sup_{a>0}M(a,\tau)$. Since $M$ is bounded from above with respect to $\tau$, uniformly in $a$, continuous with respect to $a$, such that $M(0,\tau)=0$ and $\lim_{a\to\infty}M(a,\tau)\to8\,\pi$, $M^*(\tau)$ is well defined and finite. The theorem is then a straightforward consequence of Theorem \ref{the:existence} and Proposition~\ref{lem:estM1}. \qed

\section{Numerical results}
\label{Sec:Numerics}

In this section, we numerically illustrate the above results. In particular, we show the existence of positive forward self-similar solutions with mass above $8\,\pi$ and their multiplicity when $\tau$ is large enough. We follow two different approaches: first the formulation \eqref{KSwr}--\eqref{cond3}, and then the cumulated densities formulation based on \eqref{phi2}--\eqref{S}.

\subsection{Bifurcation diagrams} The computations giving rise to Figs. 1 and 2 are based on the parametrization provided by \eqref{KSwr}--\eqref{cond3}. Numerically, one has to be careful with the origin and solve \eqref{KSwr} on the interval $(\varepsilon,\infty)$ with the initial conditions
\[
w(\varepsilon\,;s)=s-\frac14\,\varepsilon^2\,{\rm e}^s\quad \mbox{and}\quad w'(\varepsilon\,;s)=-\frac12\,\varepsilon\,{\rm e}^s\,,
\]
obtained by the Taylor expansion at $\varepsilon>0$, small enough, thus dropping higher order terms in $\varepsilon$. Observe that by \eqref{KSwr} $w''(0\,;s)=-{\rm e}^s/2$. In case of Fig. 2, one has to compute $M(s)$, which is given by \eqref{eq:M(s)}, by solving $M'(r)=2\,\pi\,{\rm e}^{w(r;s)}\,{\rm e}^{-r^2/4}\,r$ with the approximate initial condition $M(\varepsilon)=\pi\,\varepsilon^2\,{\rm e}^s$.

In Fig. 2, we recover that $M(s)\to 8\,\pi$ as $s\to \infty$. Moreover, for $\tau$ large enough, there are two solutions corresponding to a given $M$ larger than $8\,\pi$, with $M-8\,\pi$ not too large. Since it is of interest to decide for which values of $\tau$ solutions may have mass larger than $8\,\pi$, the small rectangle in Fig. 2 (left) is enlarged in Fig.~2 (right).

It can be numerically checked that solving the equations on $(\varepsilon, r_{\rm max})$ with \hbox{$r_{\rm max}\!=\!10$} gives a~good approximation of the solution. Furthermore, here we took \hbox{$\varepsilon=10^{-8}$} and $s\in[-10,20]$.
\begin{figure}[ht]\begin{center}\includegraphics[width=6cm]{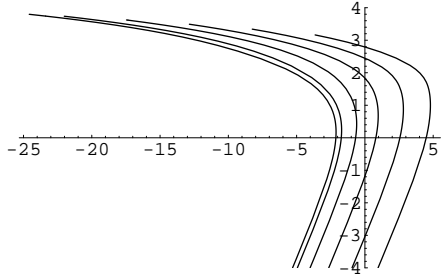}\caption{The set of all positive solutions of $\Delta v_\sigma+\frac\tau{2}\,\xi\cdot\nabla v_\sigma+\sigma\,{\rm e}^{v_\sigma}{\rm e}^{-|\xi|^2/4}=0$ in $C^2_0(\R^2)$, where $\sigma=\sigma(s)={\rm e}^{w(\infty;s)}$, is represented by the multivalued diagram $s\mapsto(\log\sigma,\log v_\sigma(0))$ for $\tau=10^{\alpha}$, $\alpha=-2$, $-1$, \ldots, $3$. Recall that the solutions $v_\sigma$ are radial and decreasing so that $v_\sigma(0)=\|v_\sigma\|_{L^\infty(\R^2)}$. We observe that $\max\limits_{s\in\R}\log\sigma(s)$ appears as an increasing function of $\tau$.}\end{center}\end{figure}

\begin{figure}[ht]\begin{center}\includegraphics[width=5.9cm]{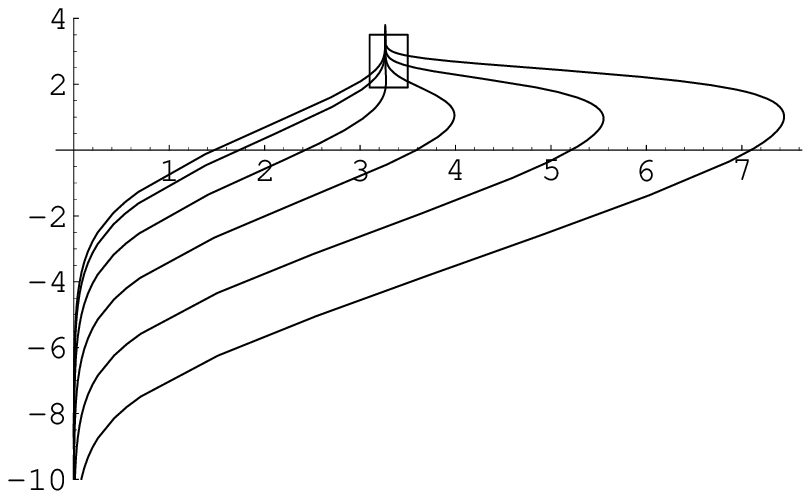}\hspace*{0.25cm}\includegraphics[width=5.9cm]{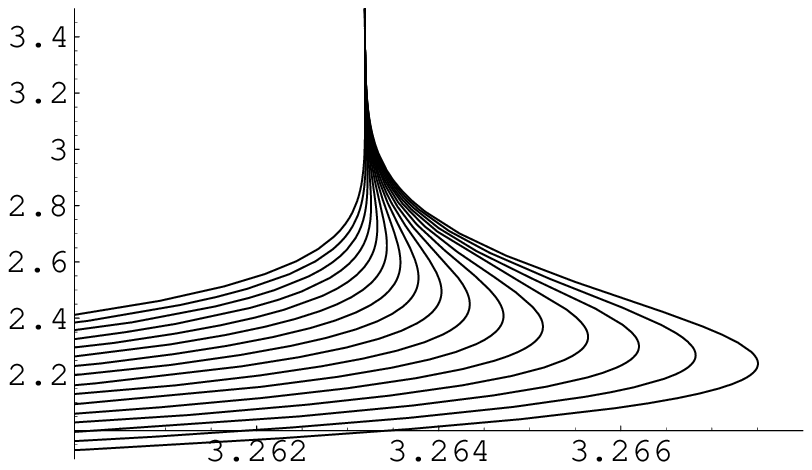}\caption{{\em Left:} The set of all positive solutions of $\Delta v_\sigma+\frac\tau{2}\,\xi\cdot\nabla v_\sigma+\sigma\,{\rm e}^{v_\sigma}{\rm e}^{-|\xi|^2/4}=0$ in $C^2_0(\R^2)$ is now represented by the diagram $s\mapsto(\log(1+M(s)),\log v_\sigma(0))$ for $\tau=10^{\alpha}$, \hbox{$\alpha=-2$, $-1$, \ldots, $3$}. We observe that $\max\limits_{s\in\R} M(s)$ appears as an increasing function of $\tau$.\newline
{\em Right:} The plot is an enlargement of the rectangle of Fig. 2 (left), with $\tau=0.60$, $0.62$, $0.64$, \ldots, $0.90$. Numerically, the first solution with mass larger than $8\,\pi$ appears for $\tau\in(0.62,0.64)$, which is far below the bound found in Section~\ref{Sec:Estimate}. This is not easy to read on the above figure, but it can be shown graphically by enlarging it enough.}\end{center}\end{figure}

\subsection{Cumulated densities} Plots and bifurcation diagrams of forward self-similar solutions can be computed in the framework of cumulated densities \eqref{phi2}--\eqref{S}, \eqref{IC}. However, again one has to be careful with the singularity at the origin. As above, since for $\eps>0$ small enough, $S'\sim\phi'\sim a$ on $(0,\eps)$ and so
\[
S(y)=a\,y+{\mathcal O}(\eps^2)\quad\mbox{and}\quad\phi''(y)\sim -\frac a4\,(1+2\,a)+{\mathcal O}(\eps)\,,
\]
we practically solve \eqref{phi2}--\eqref{S} on $(\eps,y_{\rm max})$ with the initial data
\[
\phi'(\eps)=a-\frac a4\,(1+2\,a)\,\eps\,,\quad\phi(\eps)=a\,\eps-\frac a8\,(1+2\,a)\,\eps^2\quad\mbox{and}\quad S(\eps)=a\,\eps\,.
\]
for any $y\in(0,\varepsilon)$. Obviously, having fixed $\eps>0$, one has to take $a$ in such a way that $\phi'(\eps)-a=o(a)$. Here, we choose $\varepsilon=10^{-6}$. Finally, we shall approximate $M$ from below by $\phi(y_{\rm max})$ with $y_{\rm max}$ large enough. Figs. 3 and 4 correspond to the cases $\tau=0.1$ and $\tau=10$ respectively. For $\tau=0.1$, the value $8\,\pi$ for the total mass is achieved only asymptotically in the limit $a\to\infty$. For $\tau=10$, self-similar solutions with mass $M$ larger than $8\,\pi$ exist for~$a$ large enough. Finally, Figs. 5 and 6 show the total mass as a function of $a$ and $\tau$.

\begin{figure}[hb]\begin{center}\includegraphics[height=3cm]{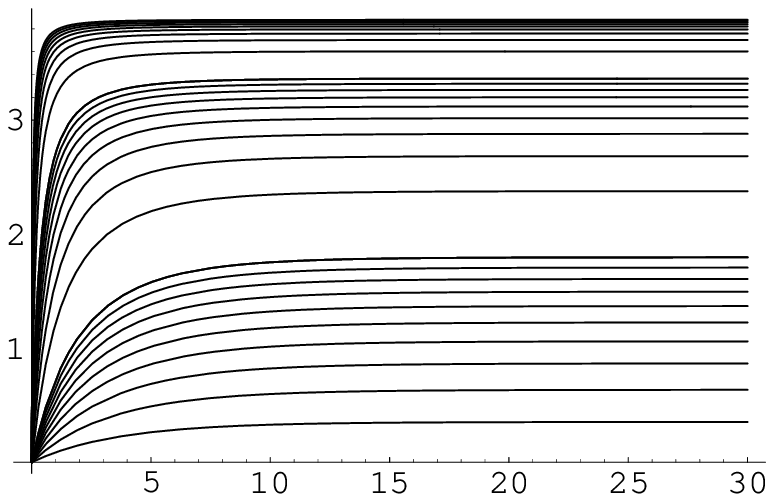}\hspace*{2cm}\includegraphics[height=3cm]{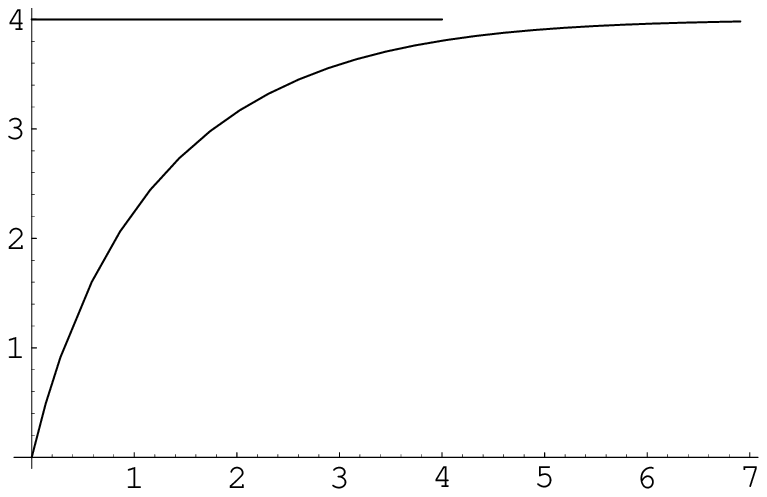}\caption{Left: Plots of $\phi$ for $\phi'(0)=a$, with $a=10^b\,c$, $b=-1$, $0$, $1$, $c\in\{1, \ldots, 10\}$ for $\tau=0.1$. Right: Plot of $b\mapsto\phi(y_{\rm max})$ in the logarithmic scale, with $\phi'(0)=a$, $a={\rm e}^b-1$, $y_{\rm max}=30$.}\end{center}\end{figure}

\begin{figure}[ht]\begin{center}\includegraphics[height=3cm]{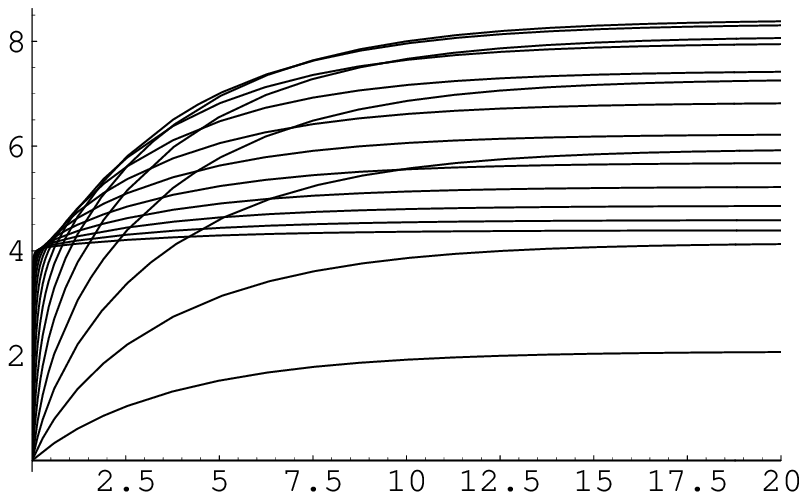}\hspace*{2cm}\includegraphics[height=3cm]{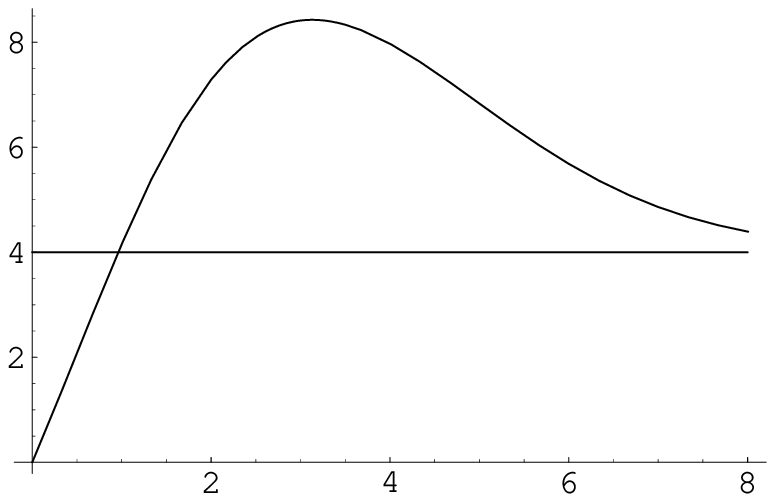}\caption{Left: Plots of $\phi$ for $\phi'(0)={\rm e}^\alpha$, with $\alpha=1$, $2$, \ldots, $20$ for $\tau=10$. Right: Plot of $\phi(y_{\rm max})$ as a function of $b$ (in the logarithmic scale), with $\phi'(0)=a$, $a={\rm e}^b-1$. Here $\tau=10$, $y_{\rm max}=30$.}\end{center}\end{figure}

\begin{figure}[ht]\begin{center}\includegraphics[height=3cm]{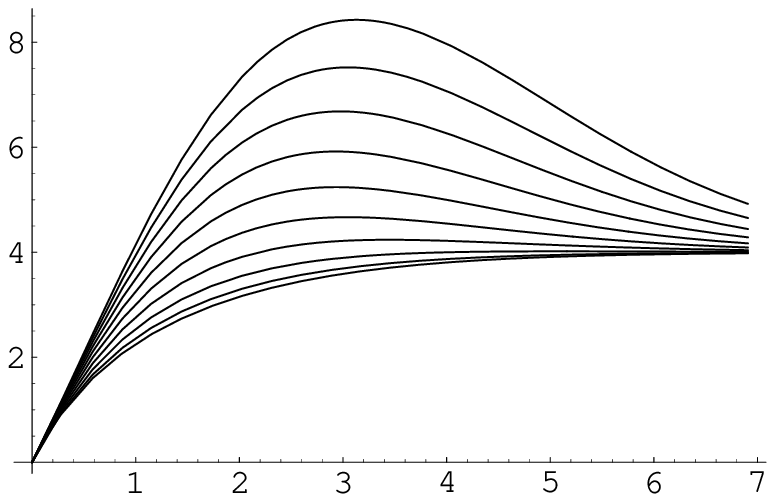}\hspace*{2cm}\includegraphics[height=3cm]{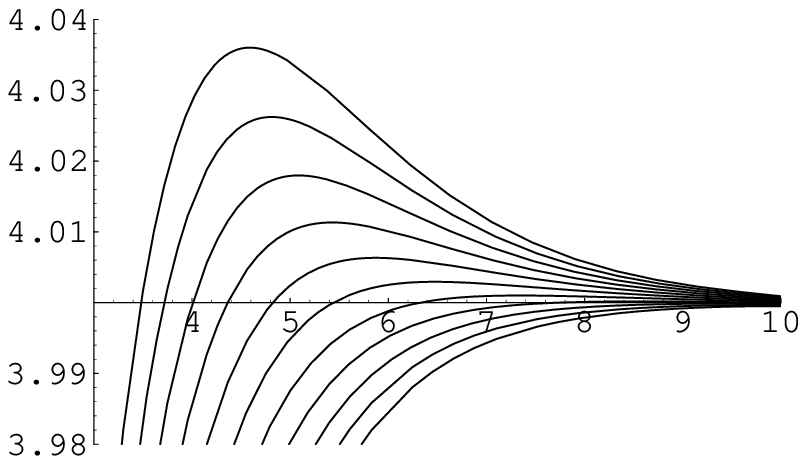}\caption{Left: The value of mass $\phi(\infty)=M(a,\tau)/(2\,\pi)$ in the logarithmic scale as a function of $a$, for $\tau=0.1\,k^2$ with $k=1$, $2$, \ldots, $10$. Right: An enlargement around the value $M(a,\tau)/(2\,\pi)=4$ in the logarithmic scale as a function of $a$, for $\tau=0.50$, $0.55$, $0.60$, \ldots, $1.00$.}\end{center}\end{figure}

\begin{figure}[ht]\begin{center}\includegraphics[height=3cm]{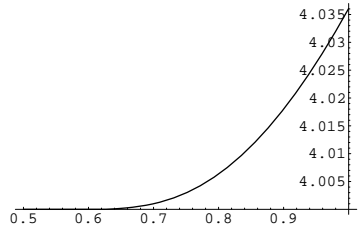}\caption{The value of the maximal (in terms of $a$) mass $\phi(\infty)=M_*(\tau)/(2\,\pi)$ as a function of~$\tau$. Numerically, the first solution with mass larger than $8\,\pi$ appears for $\tau\in(0.62,0.64)$, as already noticed at the level of Fig. 2 (right). This is again not easy to read on the above figure, but it can be shown graphically by enlarging it enough.}\end{center}\end{figure}

\section{Conclusions}
\label{Sec:Conclusions}

Self-similar solutions are much more than an example of a family of solutions. The experience of various nonlinear diffusion equations shows that they are likely to be attracting a whole class of solutions, although this is still an open question for the parabolic-parabolic Keller--Segel model with large mass (see \cite{N06} for a result for small mass solutions). It is quite reasonable to expect that well chosen perturbations of these solutions asymptotically converge in self-similar variables to the stationary solutions we have found. This actually raises a much more interesting question, which is how to determine the basin of attraction of these self-similar solutions and to understand where is the threshold between solutions for which diffusion predominates and solutions which aggregate. Clearly, it is not going to be as simple as in the parabolic-elliptic case, where a single parameter, the total mass, determines the asymptotic regime. We can conjecture that blowup occurs for mass large enough and even, maybe, as soon as the total mass of the system is above $8\,\pi$ if initial data are sufficiently concentrated.

The model considered in this paper is by many aspects ridiculously simple. See, for instance, \cite{MR2448428} to get a taste of the variety of the nonlinearities that make sense even for a~rather crude modelling purpose. Still, these models, in limiting regimes, asymptotically exhibit scaling properties similar to the ones of the parabolic-parabolic Keller--Segel model considered here. Therefore, we believe that the information gathered above, together with the methods that have been introduced, for instance, the cumulated densities reformulation of the model, should definitely be some valuable piece of information in the study of the asymptotic behaviors of the equations used in chemotaxis.

\begin{acknowledgements} The authors have been supported by Polonium contract nr.~13886SG (2007--2008). This work has been initiated during the {\sl Special semester on quantitative biology analyzed by mathematical methods}, October 1st 2007 -- January 27th, 2008, organized by RICAM, Austrian Academy of Sciences, in Linz. More recently, this research has been partially supported by the ANR CBDif-Fr, the European Commission Marie Curie Host Fellowship for the Transfer of Knowledge ``Harmonic Analysis, Nonlinear Analysis and Probability'' MTKD-CT-2004-013389, and by the Polish Ministry of Science (MNSzW) grant -- project \hbox{N201 022 32/0902}.\\
{\sl\scriptsize \copyright~2009 by the authors. This paper may be reproduced, in its entirety, for noncommercial purposes.}\end{acknowledgements}

\bibliographystyle{spmpsci}
\bibliography{References}

\end{document}